\documentclass[11pt,reqno]{amsart}
\usepackage{amsmath,amssymb,amsthm,mathrsfs}
\usepackage{geometry}
\usepackage{hyperref}
\usepackage{enumitem}
\newtheorem{theorem}{Theorem}[section]
\newtheorem{proposition}[theorem]{Proposition}
\newtheorem{lemma}[theorem]{Lemma}
\newtheorem{corollary}[theorem]{Corollary}
\theoremstyle{remark}
\newtheorem{remark}[theorem]{Remark}
\newtheorem{definition}[theorem]{Definition}

\title[Multiplicity of Positive Solutions...]
{Multiplicity of Positive Solutions for Singular Anisotropic Elliptic Equations with Variable Exponents}

\author[A.A.Esmailnejad]{Ali Akbar Esmailnejad$^1$}
\email{Hamidesmailnegad@gmail.com}
\author[S.Khademloo]{Somayeh Khademloo$^1$$^\ast$}
\email{somayehk@kth.se}
\author[S.H.Rasouli]{Seyyed Hashem Rasouli$^1$}
\email{s.h.rasouli@nit.ac.ir}
\address{$^1$Department of Basic Sciences, Babol Noshirvani University of Technology, Babol, Iran}

\thanks{$^\ast$Corresponding author}

\keywords{Anisotropic operators, Singular nonlinearities, Nehari manifold}
\subjclass[2020]{35J60, 35J75, 46E35}

\begin{document}

\begin{abstract}
We study anisotropic elliptic Dirichlet problems with variable exponents, singular nonlinearities, and subcritical superlinear reactions.
Our contribution is a direct variational construction that combines Nehari geometry and fibering maps without requiring global differentiability of the singular energy, and the result establishes the existence of two positive weak solutions for sufficiently small \(\lambda\), distinguished by their energy signs.
\end{abstract}
\maketitle

\tableofcontents

\section{Introduction}

In heterogeneous materials, the resistance to diffusion may depend
not only on the position at which a process takes place but also on
the direction in which it evolves. This feature leads to
nonlinear elliptic equations with nonstandard growth. Variable-exponent Lebesgue and Sobolev spaces provide a natural
functional framework for these models, allowing the growth of the
underlying energy to vary from point to point; see
\cite{Diening-Harjulehto-Hasto-Ruzicka,chen}.

An extension is obtained when the growth
exponent depends simultaneously on the spatial position and the
direction of differentiation. In this setting, the classical
$p(x)$-Laplacian is replaced by an anisotropic variable-exponent operator of the form
\(
\sum_{i=1}^{N}
\partial_i\left(
|\partial_i u|^{p_i(x)-2}\partial_i u
\right),
\)
where $\vec p(x)=(p_1(x),\ldots,p_N(x))$ represents the
direction-dependent growth exponents. Fundamental contributions to this
theory can be found in \cite{kov,mpr1,kha}. More recently, in
\cite{ChemsEddineRagusaRepovs2024} critical embeddings
and a concentration--compactness principle for anisotropic
variable-exponent Sobolev spaces are developed. Subsequent investigations in
\cite{ChemsEddineRagusa2025,ChemsEddineRagusa2026}
provided further compactness criteria, refined embedding results,
and applications to anisotropic elliptic equations with critical
growth.

The main difficulty in this work is that the equation contains a singular reaction term of the form $u^{-\delta(x)}$, where
$0<\delta(x)<1$, which is unbounded as $u$ approaches zero. Although such terms have a
sublinear primitive, they introduce difficulties associated with
differentiability, positivity, and the behavior of solutions near
their zero sets. When a singular contribution is coupled with a
superlinear reaction, the corresponding variational problem
exhibits competing effects at small and large amplitudes.
The singular term governs the behavior near the origin, whereas
the superlinear contribution may cause the energy functional to
be unbounded below on the ambient Sobolev space. Consequently,
standard unconstrained minimization is generally insufficient.
See \cite{10,Blasio,BaiPapageorgiouZeng2024}.

Recent investigations have considerably broadened the scope of
existence and multiplicity results for anisotropic equations.
We refer the reader to \cite{GuefaifiaBoulaarasJan2025,ZuoBahadiAkdimWu2026} . Although these works
are formulated for different operators and nonlinear structures,
they highlight the growing interest in the interaction between
anisotropy and singular effects.

Recent studies employ Nehari manifolds and fibering maps to establish multiplicity for singular multi-phase and mixed local and nonlocal variable-exponent problems (\cite{Avci2026,MalhotraPandeySreenadh2025}), while existence and uniqueness under $p$-sublinear perturbations have also been investigated in \cite{EspositoOlivaVecchi2026}.
The classical fibering approach (\cite{Drabek-Pohozaev}) and related Nehari manifold techniques (\cite{kha1,kha2}) provide an effective framework for identifying distinct critical points through the decomposition of the energy geometry.
In anisotropic variable-exponent settings, this decomposition requires refined estimates to accommodate the interplay of directional growth, singular nonlinearities, and superlinear reactions.

\medskip

Motivated by these developments, in the present paper we establish the existence of
multiple positive weak solutions under suitable assumptions
on the exponents and weight functions. In particular, we
prove that there exists a constant $\lambda_0>0$ such that,
for every $\lambda\in(0,\lambda_0)$, the problem admits
at least two distinct positive weak solutions. The central
idea is to study the fibering maps associated with the
energy functional and to decompose the Nehari manifold into
the subsets $\mathcal N_\lambda^+$,
$\mathcal N_\lambda^-$, and $\mathcal N_\lambda^0$.
For sufficiently small values of $\lambda$, the degenerate
component is excluded, allowing the two nondegenerate
components to be treated separately.

\section{Anisotropic Variable Exponent Spaces}
\label{sec:spaces}

Let us recall the properties of variable exponent spaces and anisotropic
Sobolev spaces that will be used throughout the paper. Let
$\Omega\subset\mathbb R^N$, $N\geq2$, be a bounded domain with $C^2$
boundary, and put
\[
P(\overline\Omega)
=\{p\in C(\overline\Omega):p(x)>1\text{ for every }x\in\overline\Omega\}.
\]
For any $p \in P(\bar{\Omega})$, we write
\[
p^-:=\min_{x\in\overline\Omega}p(x),\qquad
p^+:=\max_{x\in\overline\Omega}p(x).
\]

The variable exponent Lebesgue space is
\[
L^{p(x)}(\Omega)
=\left\{u\text{ measurable}:\int_\Omega |u|^{p(x)}\,dx<\infty\right\},
\]
equipped with the Luxemburg norm
\[
\|u\|_{p(x)}
=\inf\left\{\tau>0:
\int_\Omega\left|\frac{u}{\tau}\right|^{p(x)}\,dx\leq1\right\}.
\]
This space is separable and reflexive. If
$\rho_p(u)=\int_\Omega |u|^{p(x)}\,dx$, then
\begin{equation*}
\begin{aligned}
\|u\|_{p(x)}\leq1
&\ \Longrightarrow\quad
\|u\|_{p(x)}^{p^+}\leq\rho_p(u)\leq\|u\|_{p(x)}^{p^-},\\
\|u\|_{p(x)}\geq1
&\ \Longrightarrow\quad
\|u\|_{p(x)}^{p^-}\leq\rho_p(u)\leq\|u\|_{p(x)}^{p^+}.
\end{aligned}
\end{equation*}
Moreover, norm convergence to zero is equivalent to modular convergence
to zero. If $p'(x)=p(x)/(p(x)-1)$, the variable exponent H\"older inequality
gives
\[
\left|\int_\Omega uv\,dx\right|
\leq2\|u\|_{p(x)}\|v\|_{p'(x)}.
\]
We refer to \cite{FanZhao,Diening-Harjulehto-Hasto-Ruzicka} for these
properties and the corresponding Sobolev-space theory.

The space $W^{1,p(x)}(\Omega)$ consists of the functions
$u\in L^{p(x)}(\Omega)$ such that $|\nabla u|\in L^{p(x)}(\Omega)$.
It may be equipped with the norm
$\|u\|_{p(x)}+\|\nabla u\|_{p(x)}$.

\subsection{Preliminaries}

Let $\vec p(x)=(p_1(x),\ldots,p_N(x))$, where
$p_i\in P(\overline\Omega)$, and define
\[
p^-:=\min_{1\leq i\leq N}p_i^-,\qquad
p^+:=\max_{1\leq i\leq N}p_i^+,
\]
\[
p_m(x):=\min_{1\leq i\leq N}p_i(x),\qquad
p_M(x):=\max_{1\leq i\leq N}p_i(x).
\]
The modular associated with our setting is
\begin{equation*}
\rho(u):=\sum_{i=1}^N\int_\Omega
\bigl(|\partial_i u|^{p_i(x)}+|u|^{p_i(x)}\bigr)\,dx
\end{equation*}
and its Luxemburg norm
\begin{equation*}
\|u\|_X:=\inf\{\tau>0:\rho(u/\tau)\leq1\}.
\end{equation*}
Throughout the paper,
\[
X=W_0^{1,\vec p(x)}(\Omega)
:=\overline{C_c^\infty(\Omega)}^{\|\cdot\|_X}.
\]
Thus, the zero boundary condition is understood through this closure.
The space $X$ is a reflexive Banach space, and
$X\hookrightarrow W_0^{1,p^-}(\Omega)$ continuously. It is stable under
absolute values, with
\[
\rho(|u|)=\rho(u),\qquad \||u|\|_X=\|u\|_X.
\]
These facts follow by viewing $X$ as a closed subspace of the finite
product of the variable exponent spaces containing $u$ and its
coordinate derivatives. The Sobolev chain rule gives the assertion
concerning absolute values; see also \cite{FanAnisotropic2011}.

\begin{proposition}\label{prop:modular}
For every $u\in X$, the following inequalities hold:
\begin{equation}\label{eq:Xmodular}
\begin{aligned}
\|u\|_X\leq1
&\ \Longrightarrow\quad
\|u\|_X^{p^+}\leq\rho(u)\leq\|u\|_X^{p^-},\\
\|u\|_X\geq1
&\ \Longrightarrow\quad
\|u\|_X^{p^-}\leq\rho(u)\leq\|u\|_X^{p^+}.
\end{aligned}
\end{equation}
In particular, $\rho(u)=1$ if and only if $\|u\|_X=1$.
\end{proposition}

\begin{proof}
For $u\neq0$, the definition of the Luxemburg norm gives
$\rho(u/\|u\|_X)=1$. Factoring $\|u\|_X^{p_i(x)}$ from each integrand
and using $p^-\leq p_i(x)\leq p^+$ proves \eqref{eq:Xmodular}.
The case $u=0$ is immediate.
\end{proof}

Define
\[
p_m^*(x):=
\begin{cases}
\displaystyle\frac{Np_m(x)}{N-p_m(x)},&p_m(x)<N,\\[2mm]
+\infty,&p_m(x)\geq N.
\end{cases}
\]
Only strictly subcritical embeddings are needed below.

\begin{proposition}
\label{prop:embedding}
If $r\in C(\overline\Omega)$ and
\[
1\leq r(x)<p_m^*(x)\qquad\text{for every }x\in\overline\Omega,
\]
then $X\hookrightarrow L^{r(x)}(\Omega)$ continuously and compactly.
\end{proposition}

\subsection{Estimates for the nonlinear terms}

Let $\delta,q\in C(\overline\Omega)$ satisfy
\[
0<\delta^-\leq\delta(x)\leq\delta^+<1,
\qquad 1<q(x)<p_m^*(x),
\]
and let $f,g\in L^\infty(\Omega)$ be nonnegative. It is convenient to put
\begin{equation*}
s(x):=1-\delta(x),\qquad
s^-:=1-\delta^+,\qquad s^+:=1-\delta^-.
\end{equation*}
Thus $0<s^-\leq s(x)\leq s^+<1$. Define
\begin{align*}
F(u)&:=\int_\Omega f(x)|u|^{s(x)}\,dx,
& S(u)&:=\int_\Omega\frac{f(x)}{s(x)}|u|^{s(x)}\,dx,\\
G(u)&:=\int_\Omega g(x)|u|^{q(x)}\,dx,
& Q(u)&:=\int_\Omega\frac{g(x)}{q(x)}|u|^{q(x)}\,dx.
\end{align*}
We work directly with these integrals. 

\begin{proposition}\label{prop:pointwise-estimates}
There exist constants $C_f,C_g>0$ such that
\begin{equation}\label{eq:pointwise-delta}
F(u)\leq C_f
\begin{cases}
\|u\|_X^{s^-},&\|u\|_X\leq1,\\
\|u\|_X^{s^+},&\|u\|_X\geq1,
\end{cases}
\end{equation}
and
\begin{equation}\label{eq:pointwise-q}
G(u)\leq C_g
\begin{cases}
\|u\|_X^{q^-},&\|u\|_X\leq1,\\
\|u\|_X^{q^+},&\|u\|_X\geq1.
\end{cases}
\end{equation}
The functionals $F,S,G,Q$ are sequentially weakly continuous on $X$.
Moreover,
\begin{equation}\label{eq:primitive-bounds}
\frac{F(u)}{s^+}\leq S(u)\leq\frac{F(u)}{s^-},
\qquad
\frac{G(u)}{q^+}\leq Q(u)\leq\frac{G(u)}{q^-}.
\end{equation}
\end{proposition}

\begin{proof}
The embedding into $L^1(\Omega)$ and the inequality
$|v|^{s(x)}\leq1+|v|$ give a uniform bound for $F(v)$ on
$\|v\|_X=1$. The embedding into $L^{q(x)}(\Omega)$ gives the analogous
bound for $G(v)$. Writing $u=\|u\|_Xv$ proves
\eqref{eq:pointwise-delta} and \eqref{eq:pointwise-q}.

If $u_n\rightharpoonup u$ in $X$, Proposition~\ref{prop:embedding}
gives strong convergence in $L^1(\Omega)$ and $L^{q(x)}(\Omega)$.
For $0<s<1$,
\[
\bigl||a|^s-|b|^s\bigr|\leq|a-b|^s.
\]
Consequently,
\[
|F(u_n)-F(u)|
\leq\|f\|_\infty\int_\Omega
\bigl(|u_n-u|^{s^-}+|u_n-u|^{s^+}\bigr)\,dx\longrightarrow0,
\]
where the last assertion follows from strong $L^1$ convergence and
H\"older's inequality with constant exponents. The same argument applies
to $S$, since $s^->0$. Strong $L^{q(x)}$ convergence implies convergence
of the $q(x)$-power integrals, and hence of $G$ and $Q$.
Finally, \eqref{eq:primitive-bounds} follows from the pointwise bounds
on $s(x)$ and $q(x)$.
\end{proof}

\section{Statement of the Singular Problem and Variational Setting}
\label{sec:problem}

In this section, we formulate the singular anisotropic elliptic problem
that will be analyzed throughout the paper. For a given parameter
$\lambda>0$, we consider
\begin{equation}\label{P}
\begin{cases}
\displaystyle
-\sum_{i=1}^N\partial_i\bigl(|\partial_i u|^{p_i(x)-2}\partial_i u\bigr)
+\sum_{i=1}^N|u|^{p_i(x)-2}u
=\lambda f(x)u^{-\delta(x)}+g(x)u^{q(x)-1}
&\text{in }\Omega,\\[2mm]
u>0&\text{a.e. in }\Omega,\\
u=0&\text{on }\partial\Omega.
\end{cases}
\end{equation}
Here $\Omega\subset\mathbb R^N$, $N\geq2$, is a bounded domain with
$C^2$ boundary. We impose the following assumptions:
\begin{align*}
\mathrm{(A1)}\quad&
p_i,\delta,q\in C(\overline\Omega),\qquad
1<p^-\leq p^+<\infty,\qquad
0<\delta^-\leq\delta^+<1;\\
\mathrm{(A2)}\quad&
0<1-\delta^+\leq1-\delta^-<p^-\leq p^+<q^-\leq q^+,\\
& q(x)<p_m^*(x)\qquad\text{for every }x\in\overline\Omega;\\
\mathrm{(A3)}\quad&
f,g\in C(\overline\Omega),\qquad
f>0\text{ a.e. in }\Omega,\qquad
g\geq0,\quad g\not\equiv0.
\end{align*}

\begin{definition}\label{def:weak}
A function $u\in X$ is a weak solution of \eqref{P} if $u>0$ a.e. in
$\Omega$,
\begin{equation}\label{eq:singular-integrability}
\int_\Omega f(x)u^{-\delta(x)}|\varphi|\,dx<\infty
\qquad\text{for every }\varphi\in X,
\end{equation}
and
\begin{equation}\label{weak}
\begin{aligned}
&\sum_{i=1}^N\int_\Omega
|\partial_i u|^{p_i(x)-2}\partial_i u\,\partial_i\varphi\,dx
+\sum_{i=1}^N\int_\Omega u^{p_i(x)-1}\varphi\,dx\\
&\hspace{12mm}=
\lambda\int_\Omega f(x)u^{-\delta(x)}\varphi\,dx
+\int_\Omega g(x)u^{q(x)-1}\varphi\,dx
\qquad\text{for every }\varphi\in X.
\end{aligned}
\end{equation}
\end{definition}

The energy functional associated with \eqref{P} is
\begin{equation*}
\mathcal I_\lambda(u)=E(u)-\lambda S(u)-Q(u),\qquad u\in X,
\end{equation*}
where
\begin{equation*}
E(u):=\sum_{i=1}^N\int_\Omega
\frac{|\partial_i u|^{p_i(x)}+|u|^{p_i(x)}}{p_i(x)}\,dx.
\end{equation*}
Note that the functionals $E$ and $Q$ belong to
$C^1(X,\mathbb R)$, and $E$ is convex. In particular,
\begin{equation*}
\langle E'(u),\varphi\rangle
={}\sum_{i=1}^N\int_\Omega
|\partial_i u|^{p_i(x)-2}\partial_i u\,\partial_i\varphi\,dx
+\sum_{i=1}^N\int_\Omega |u|^{p_i(x)-2}u\varphi\,dx,
\end{equation*}
and $\langle E'(u),u\rangle=\rho(u)$. The functional
$\mathcal I_\lambda$ is continuous and sequentially weakly lower
semicontinuous, by the convexity of $E$ and
Proposition~\ref{prop:pointwise-estimates}.

The superlinear term prevents global coercivity. Indeed, choosing
$\psi\in C_c^\infty(\Omega)$ with $G(\psi)>0$, we have
$\mathcal I_\lambda(t\psi)\to-\infty$ as $t\to\infty$ because
$q^->p^+$. Note that $S$ is not differentiable on $X$. Therefore, we do not assume
that $\mathcal I_\lambda\in C^1(X,\mathbb R)$ and do not apply a
Lagrange multiplier argument to it.

Instead, we use the differentiability of the restrictions
$t\mapsto\mathcal I_\lambda(tu)$ for $t>0$, construct minimizers on
the two Nehari components, and identify their weak equations by
one-sided variations. This direct treatment of the singular energy
is consistent with the constrained variational approach used in
\cite{kamel,PapageorgiouRepovsVetro2021}, while the arguments required
for the present anisotropic variable-exponent setting are given below.
The construction is carried out directly for the singular functional.

Here is the main result of this paper.

\begin{theorem}\label{thm:two-solutions}
Assume {\rm(A1)--(A3)}. There exists $\lambda_0>0$ such that, for every
$\lambda\in(0,\lambda_0)$, problem \eqref{P} has at least two distinct
positive weak solutions $u_\lambda,v_\lambda\in X$. Moreover,
\[
\mathcal I_\lambda(u_\lambda)<0<\mathcal I_\lambda(v_\lambda).
\]
Positivity is understood almost everywhere, as in
Definition~\ref{def:weak}.
\end{theorem}

The constant $\lambda_0$ will be chosen explicitly in terms of the
embedding estimates in Section~\ref{sec:Nehari}. The proof of the theorem
is completed in Section~\ref{sec:main-proof}.

\section{Nehari Manifold and Fibering Maps}
\label{sec:Nehari}

For $u\in X\setminus\{0\}$, define the fibering map
\[
\Phi_u(t):=\mathcal I_\lambda(tu),\qquad t\geq0.
\]
Although $\mathcal I_\lambda$ is not assumed to be differentiable on
$X$, the function $\Phi_u$ belongs to $C^2((0,\infty))$. To simplify
notation, put
\[
a_i(u,x):=|\partial_i u|^{p_i(x)}+|u|^{p_i(x)},\qquad
f_u(x):=f(x)|u|^{s(x)},\qquad g_u(x):=g(x)|u|^{q(x)}.
\]
Then
\begin{equation*}
\Phi_u'(t)
=\sum_{i=1}^N\int_\Omega t^{p_i(x)-1}a_i(u,x)\,dx
-\lambda\int_\Omega t^{s(x)-1}f_u\,dx
-\int_\Omega t^{q(x)-1}g_u\,dx,
\end{equation*}
and
\[\Phi_u''(t)
=\sum_{i=1}^N\int_\Omega
(p_i(x)-1)t^{p_i(x)-2}a_i(u,x)\,dx
+\lambda\int_\Omega\delta(x)t^{-\delta(x)-1}f_u\,dx
-\int_\Omega(q(x)-1)t^{q(x)-2}g_u\,dx.\]
The Nehari set is defined by the radial constraint
\begin{equation*}
\mathcal N_\lambda
:=\{u\in X\setminus\{0\}:J_\lambda(u)=0\},
\qquad
J_\lambda(u):=\Phi_u'(1)=\rho(u)-\lambda F(u)-G(u).
\end{equation*}
Thus $t\Phi_u'(t)=J_\lambda(tu)$. We retain the term
\emph{Nehari manifold}, but no global differentiable-manifold structure
is required.  Every nontrivial weak solution of \eqref{P} belongs to
$\mathcal N_\lambda$, by testing its equation with $u$.

For $u\in\mathcal N_\lambda$, set
\begin{equation*}
D_\lambda(u)
:=\sum_{i=1}^N\int_\Omega p_i(x)a_i(u,x)\,dx
-\lambda\int_\Omega s(x)f_u\,dx
-\int_\Omega q(x)g_u\,dx.
\end{equation*}
Since $J_\lambda(u)=0$, we have $D_\lambda(u)=\Phi_u''(1)$.
Accordingly,
\[
\mathcal N_\lambda^\pm
:=\{u\in\mathcal N_\lambda:D_\lambda(u)\gtrless0\},
\qquad
\mathcal N_\lambda^0
:=\{u\in\mathcal N_\lambda:D_\lambda(u)=0\}.
\]

\begin{lemma}\label{lem:coercive}
The functional $\mathcal I_\lambda$ is coercive and bounded below on
$\mathcal N_\lambda$.
\end{lemma}

\begin{proof}
For $u\in\mathcal N_\lambda$, we have
$G(u)=\rho(u)-\lambda F(u)$. Therefore,
\begin{equation*}
\begin{aligned}
\mathcal I_\lambda(u)
&\geq\frac{\rho(u)}{p^+}-\frac{\lambda F(u)}{s^-}
-\frac{G(u)}{q^-}\\
&=\left(\frac1{p^+}-\frac1{q^-}\right)\rho(u)
-\lambda\left(\frac1{s^-}-\frac1{q^-}\right)F(u).
\end{aligned}
\end{equation*}
If $\|u\|_X\geq1$, Propositions~\ref{prop:modular} and
\ref{prop:pointwise-estimates} give
\[
\mathcal I_\lambda(u)
\geq\left(\frac1{p^+}-\frac1{q^-}\right)\|u\|_X^{p^-}
-\lambda C_f\left(\frac1{s^-}-\frac1{q^-}\right)\|u\|_X^{s^+}.
\]
Since $q^->p^+$ and $p^->s^+$, RHS tends to $+\infty$
as $\|u\|_X\to\infty$. For $\|u\|_X\leq1$, the same estimates
provide a uniform lower bound. 
\end{proof}

\medskip

The following lemma shows that, for sufficiently small $\lambda>0$,
both the energy functional and the Nehari constraint functional remain strictly positive
on a suitable sphere in $X$.
This property plays an important role in the analysis of fibering maps
and the construction of distinct variational solutions.

\begin{lemma}\label{lem:sphere}
Choose
\begin{equation}\label{eq:threshold}
r_0:=\min\left\{\frac12,
\left(\frac1{4C_g}\right)^{1/(q^--p^+)}\right\},\qquad
\lambda_0:=\frac{s^-}{4p^+C_f}\,r_0^{p^+-s^-},
\qquad
\eta_0:=\frac{r_0^{p^+}}{2p^+}.
\end{equation}
For every $\lambda\in(0,\lambda_0)$ and every $u\in X$ with
$\|u\|_X=r_0$, one has
\begin{equation*}
J_\lambda(u)>0,\qquad \mathcal I_\lambda(u)\geq\eta_0>0.
\end{equation*}
\end{lemma}

\begin{proof}
Since $r_0<1$, the preceding estimates give
\[
\rho(u)\geq r_0^{p^+},\qquad
F(u)\leq C_fr_0^{s^-},\qquad G(u)\leq C_gr_0^{q^-}.
\]
By the choice of $r_0$ and $\lambda_0$,
\[
C_gr_0^{q^-}\leq\frac14r_0^{p^+},\qquad
\lambda C_fr_0^{s^-}<\frac{s^-}{4p^+}r_0^{p^+}.
\]
It follows that
\[
J_\lambda(u)>
\left(1-\frac14-\frac{s^-}{4p^+}\right)r_0^{p^+}>0.
\]
Similarly,
\[
\mathcal I_\lambda(u)
\geq\left(\frac1{p^+}-\frac1{4p^+}-\frac1{4q^-}\right)r_0^{p^+}
\geq\frac{r_0^{p^+}}{2p^+}.
\]
\end{proof}

\begin{lemma}[Geometry of the fibering maps]\label{lem:fibering}
Let $\lambda\in(0,\lambda_0)$ and $u\in X\setminus\{0\}$.
If $G(u)>0$, there exist exactly two numbers
\[
0<t^+(u)<\frac{r_0}{\|u\|_X}<t^-(u)
\]
such that
$t^+(u)u\in\mathcal N_\lambda^+$ and
$t^-(u)u\in\mathcal N_\lambda^-$.
If $G(u)=0$, there is exactly one such number $t^+(u)$, and
$t^+(u)u\in\mathcal N_\lambda^+$.
All these radial intersections are nondegenerate.
\end{lemma}

\begin{proof}
Assumption {\rm(A3)} implies $F(u)>0$ for every $u\neq0$.
Set
\(K_u(t):=t^{-s^+}J_\lambda (tu),\) for all \(t>0.\)
Differentiation gives
\begin{equation}\label{eq:L}
\begin{aligned}
K_u'(t)
&=t^{p^+-s^+-1}
\Big(\sum_{i=1}^N\int_\Omega
(p_i(x)-s^+)a_i(u,x)t^{p_i(x)-p^+}\,dx\\
&+\lambda\int_\Omega(s^+-s(x))f_u t^{s(x)-p^+}\,dx
-\int_\Omega(q(x)-s^+)g_u t^{q(x)-p^+}\,dx\Big).
\end{aligned}
\end{equation}
Suppose first that $G(u)>0$. Since
\[
s(x)\le s^+<p^-\le p_i(x)\le p^+<q^-\le q(x),
\]
the first two terms inside the parentheses in
\eqref{eq:L} are nonincreasing in $t$, whereas
the last term is strictly decreasing.
Moreover, the expression inside the parentheses
is positive for sufficiently small $t$ and tends
to $-\infty$ as $t\to\infty$.
Consequently, it vanishes at exactly one point,
and hence $K_u$ is strictly increasing up to
that point and strictly decreasing thereafter.

For $0<t\leq1$, we have
\[
K_u(t)\leq t^{p^--s^+}\rho(u)-\lambda F(u).
\]
Since $p^->s^+$ and $F(u)>0$, it follows that
$K_u(t)<0$ for all sufficiently small $t$.
Moreover, $K_u(t)\to-\infty$ as $t\to\infty$,
because $q^->p^+$.

Lemma~\ref{lem:sphere} gives
\[
K_u\left(\frac{r_0}{\|u\|_X}\right)>0.
\]
Thus $K_u$, and hence $J_\lambda(tu)$ (as a function of $t$), has exactly two zeros, one on each side
of $r_0/\|u\|_X$. At the first zero $(J_\lambda(tu))'>0$, and at the second
$(J_\lambda(tu))'<0$. Since $\Phi_u'=J_\lambda(tu)/t$, these zeros belong to the stated
Nehari components and are nondegenerate.

If $G(u)=0$, the last term in \eqref{eq:L} vanishes, so
$K_u'(t)>0$ for all $t>0$. The function $K_u$ is negative near zero
and tends to $+\infty$. It therefore has exactly one zero, again
before $r_0/\|u\|_X$, and this zero is nondegenerate.
\end{proof}

\begin{corollary}\label{cor:geometry}
For every $\lambda\in(0,\lambda_0)$,
\[
\mathcal N_\lambda^0=\emptyset,\qquad
\mathcal N_\lambda^+\neq\emptyset,\qquad
\mathcal N_\lambda^-\neq\emptyset.
\]
Furthermore,
\begin{equation}\label{eq:branch-norms}
\|u\|_X<r_0\quad(u\in\mathcal N_\lambda^+),\qquad
\|v\|_X>r_0\quad(v\in\mathcal N_\lambda^-),
\end{equation}
and
\begin{equation}\label{eq:energy-signs}
\mathcal I_\lambda(u)<0\quad(u\in\mathcal N_\lambda^+),\qquad
\mathcal I_\lambda(v)\geq\eta_0\quad(v\in\mathcal N_\lambda^-).
\end{equation}
For every $v\in\mathcal N_\lambda^-$,
\begin{equation}\label{eq:global-ray-max}
\mathcal I_\lambda(v)=\max_{t\geq0}\mathcal I_\lambda(tv).
\end{equation}
\end{corollary}

\begin{proof}
Nondegeneracy and the norm separation follow from
Lemma~\ref{lem:fibering}. Every nonzero direction has a positive-branch
intersection. Since $g\not\equiv0$, a smooth function can be chosen
with $G(u)>0$, and its ray also meets the negative branch.

The fibering map decreases from $\Phi_u(0)=0$ to its first critical
point; hence the energy on $\mathcal N_\lambda^+$ is negative.
When $G(u)>0$, it then increases to its second critical point and
subsequently decreases to $-\infty$. The second critical value is at
least the value at $r_0u/\|u\|_X$, which is bounded below by
$\eta_0$. Since this value is positive, it is also the global maximum
on the ray, including $t=0$.
\end{proof}

\begin{lemma}\label{lem:ray-continuity}
For fixed $\lambda\in(0,\lambda_0)$, the map
$u\mapsto t^-(u)$ is continuous on
\[
\mathcal D:=\{u\in X:G(u)>0\}.
\]
\end{lemma}

\begin{proof}
The function $(t,u)\mapsto J_\lambda(tu)$ is continuous for $t>0$,
and its radial derivative is continuous there. At $t=t^-(u)$, this
radial derivative is strictly negative. Choose two nearby values of
$t$, one on each side of the zero, at which $J_\lambda(tu)$ has
opposite signs. For all directions sufficiently close to $u$, the
same signs persist. The unique descending zero for each such
direction lies between these two values. Letting them approach
$t^-(u)$ proves continuity. This argument uses only radial
differentiability, not differentiability with respect to $u$.
\end{proof}

\section{Existence of Nontrivial Minimizers}
\label{sec:minimizers}

This section presents the constrained minimization arguments on the
two nondegenerate components of the Nehari set.

\begin{theorem}\label{th5.2}
For every $\lambda\in(0,\lambda_0)$, there exists a nonnegative
$u_\lambda\in\mathcal N_\lambda^+$ such that
\[
\mathcal I_\lambda(u_\lambda)
=\inf_{u\in\mathcal N_\lambda^+}\mathcal I_\lambda(u)<0.
\]
Moreover, $u_\lambda$ is a local minimizer of $\mathcal I_\lambda$
on $X$.
\end{theorem}

\begin{proof}
Consider the closed ball
\[
\overline B_{r_0}:=\{u\in X:\|u\|_X\leq r_0\}.
\]
It is weakly compact, and $\mathcal I_\lambda$ is sequentially weakly
lower semicontinuous. The functional therefore attains its minimum
on this ball at some $u_\lambda$.

Choose a nonzero, nonnegative $\psi\in C_c^\infty(\Omega)$.
Then $F(\psi)>0$. For $0<t<1$,
\[
\mathcal I_\lambda(t\psi)
\leq\frac{\rho(\psi)}{p^-}t^{p^-}
-\frac{\lambda F(\psi)}{s^+}t^{s^+}<0
\]
for all sufficiently small $t$, since $s^+<p^-$. Hence
$\mathcal I_\lambda(u_\lambda)<0$ and $u_\lambda\neq0$.
On the boundary of the ball the energy is positive by
Lemma~\ref{lem:sphere}, so $\|u_\lambda\|_X<r_0$.
Thus $u_\lambda$ is a local minimizer on $X$.

The one-dimensional map $t\mapsto\mathcal I_\lambda(tu_\lambda)$
has a local minimum at $t=1$. In particular,
$J_\lambda(u_\lambda)=0$. By \eqref{eq:branch-norms} and the absence
of a degenerate component, $u_\lambda\in\mathcal N_\lambda^+$.
Every element of $\mathcal N_\lambda^+$ belongs to
$\overline B_{r_0}$, so the minimizer also attains the infimum on
$\mathcal N_\lambda^+$. Finally, replacing $u_\lambda$ by
$|u_\lambda|$ preserves the norm, the energy, and the Nehari component.
We may therefore take $u_\lambda\geq0$.
\end{proof}

\begin{lemma}\label{lem:negative-nonvanishing}
For every $v\in\mathcal N_\lambda^-$,
\begin{equation*}
G(v)>\frac{p^--s^+}{q^+-s^+}\rho(v)
\geq\frac{p^--s^+}{q^+-s^+}r_0^{p^+}.
\end{equation*}
\end{lemma}

\begin{proof}
Using $\lambda F(v)=\rho(v)-G(v)$, we obtain
\[
D_\lambda(v)
\geq p^-\rho(v)-\lambda s^+F(v)-q^+G(v)
=(p^--s^+)\rho(v)-(q^+-s^+)G(v).
\]
Since $D_\lambda(v)<0$, the first inequality follows. The second
follows from $\|v\|_X>r_0$, $r_0<1$, and \eqref{eq:Xmodular},
regardless of whether $\|v\|_X$ is smaller or larger than one.
\end{proof}

\begin{theorem}\label{th5.5}
For every $\lambda\in(0,\lambda_0)$, there exists a nonnegative
$v_\lambda\in\mathcal N_\lambda^-$ such that
\[
\mathcal I_\lambda(v_\lambda)
=\inf_{v\in\mathcal N_\lambda^-}\mathcal I_\lambda(v)
\geq\eta_0>0.
\]
\end{theorem}

\begin{proof}
Put
\[
c_\lambda^-:=\inf_{v\in\mathcal N_\lambda^-}\mathcal I_\lambda(v).
\]
This number is finite because the set is nonempty, and it is at least
$\eta_0$ by \eqref{eq:energy-signs}. Let
$\{v_n\}\subset\mathcal N_\lambda^-$ be a minimizing sequence.
Replacing its elements by their absolute values, we assume $v_n\geq0$.
Lemma~\ref{lem:coercive} shows that the sequence is bounded in $X$.
After passing to a subsequence,
\[
v_n\rightharpoonup v\text{ in }X,\qquad
v_n\to v\text{ in }L^{q(x)}(\Omega)\text{ and }L^1(\Omega),
\]
with $v\geq0$. Proposition~\ref{prop:pointwise-estimates} and
Lemma~\ref{lem:negative-nonvanishing} imply
\[
G(v)=\lim_{n\to\infty}G(v_n)
\geq\frac{p^--s^+}{q^+-s^+}r_0^{p^+}>0.
\]
In particular, $v\neq0$ and the descending radial parameter
$t_*=t^-(v)$ is well defined.

By weak lower semicontinuity and the global maximum property on
each ray,
\[
c_\lambda^-
\leq\mathcal I_\lambda(t_*v)
\leq\liminf_{n\to\infty}\mathcal I_\lambda(t_*v_n)
\leq\limsup_{n\to\infty}\mathcal I_\lambda(t_*v_n)
\leq\lim_{n\to\infty}\mathcal I_\lambda(v_n)
=c_\lambda^-.
\]
Thus $v_\lambda:=t_*v\in\mathcal N_\lambda^-$ attains the infimum.
\end{proof}

\section{Proof of the Main Result}
\label{sec:main-proof}

In this section, we show that the minimizers constructed above are
positive weak solutions. The essential point is to establish the
integrability of the singular term before differentiating it in an
arbitrary direction.

\begin{lemma}\label{lem:identification}
Assume  {\rm(A1)--(A3)} and let $\lambda>0$, $w\in X\setminus\{0\}$ with $w\geq0$. Suppose that for every
$\varphi\in X$ and all sufficiently small $\varepsilon>0$ there is
$t_\varepsilon>0$, possibly depending on $\varphi$, such that
$t_\varepsilon\to1$ and
\begin{equation}\label{eq:variation-property}
\mathcal I_\lambda\bigl(t_\varepsilon(w+\varepsilon\varphi)\bigr)
\geq\mathcal I_\lambda(t_\varepsilon w).
\end{equation}
Then $w$ is a positive weak solution of \eqref{P}.
\end{lemma}

\begin{proof}
Take $h\in X$, $h\geq0$.
By \eqref{eq:variation-property},
\begin{equation}\label{eq:positive-variation}
\lambda\int_\Omega
\frac{f(x)t_\varepsilon^{s(x)}}{s(x)}
\frac{(w+\varepsilon h)^{s(x)}-w^{s(x)}}{\varepsilon}\,dx
\leq
\frac{(E-Q)(t_\varepsilon(w+\varepsilon h))-(E-Q)(t_\varepsilon w)}
{\varepsilon}
\end{equation}
Since $t_\varepsilon\to1$ and $E-Q\in C^1(X)$,
the mean value formula gives
\[
\begin{aligned}
\frac{(E-Q)(t_\varepsilon(w+\varepsilon h))
 -(E-Q)(t_\varepsilon w)}{\varepsilon} =&
\int_0^1
\left\langle
(E-Q)'(t_\varepsilon w+
\theta\varepsilon t_\varepsilon h),
t_\varepsilon h
\right\rangle\,d\theta\\
\longrightarrow &
\langle(E-Q)'(w),h\rangle.
\end{aligned}
\]
The integrand on the left is nonnegative. At any point where $w=0$,
$f>0$, and $h>0$, it tends to $+\infty$, because $s(x)<1$.
Fatou's lemma therefore implies that $\{w=0,f>0,h>0\}$ has measure
zero. Taking nonnegative smooth test functions which are positive
on the members of a countable covering of $\Omega$, and using
$f>0$ a.e., we obtain $w>0$ a.e. in $\Omega$.

Another application of Fatou's lemma to
\eqref{eq:positive-variation} gives
\begin{equation}\label{eq:singular-bound}
\lambda\int_\Omega f(x)w^{-\delta(x)}h\,dx
\leq\langle(E-Q)'(w),h\rangle<\infty
\qquad(h\in X,\ h\geq0).
\end{equation}
In particular, choosing $h=|\varphi|$ proves
\eqref{eq:singular-integrability} for every $\varphi\in X$.
It also shows that the singular density defines a bounded linear
functional on $X$.

We may now return to arbitrary signed directions. For $a>0$,
$b\in\mathbb R$, and $0<s<1$, the elementary inequality
\begin{equation}\label{eq:scalar-domination}
\bigl||a+b|^s-a^s\bigr|\leq a^{s-1}|b|
\end{equation}
holds. Indeed, after dividing by $a^s$, it follows by considering
whether $|1+b/a|$ is smaller or larger than one and using the
concavity of $t\mapsto t^s$. Applying
\eqref{eq:scalar-domination} with $a=w(x)$ and
$b=\varepsilon\varphi(x)$ gives
\[
\frac{f(x)}{s(x)}
\frac{\bigl||w+\varepsilon\varphi|^{s(x)}-w^{s(x)}\bigr|}
{\varepsilon}
\leq\frac1{s^-}f(x)w^{-\delta(x)}|\varphi|.
\]
The right-hand side is integrable by \eqref{eq:singular-bound}.
The factors $t_\varepsilon^{s(x)}$ are uniformly bounded and converge
to one. Dominated convergence consequently yields
\begin{equation*}
\lim_{\varepsilon\downarrow0}
\frac{S(t_\varepsilon(w+\varepsilon\varphi))-S(t_\varepsilon w)}
{\varepsilon}=\int_\Omega f(x)w^{-\delta(x)}\varphi\,dx.
\end{equation*}
Dividing \eqref{eq:variation-property} by $\varepsilon$ after subtracting
the right-hand side, and passing to the limit, one obtains
\[
\langle E'(w),\varphi\rangle
-\int_\Omega g(x)w^{q(x)-1}\varphi\,dx
-\lambda\int_\Omega f(x)w^{-\delta(x)}\varphi\,dx\geq0.
\]
Applying the same argument to $-\varphi$ gives the reverse inequality.
Thus equality holds for every $\varphi\in X$, which is precisely
\eqref{weak}. 
\end{proof}

\medskip

\begin{proof}[Proof of Theorem~\ref{thm:two-solutions}]
Let $0<\lambda<\lambda_0$, with $\lambda_0$ given by
\eqref{eq:threshold}. By Theorem~\ref{th5.2}, there is a nonnegative
local minimizer $u_\lambda\in\mathcal N_\lambda^+$.
For any $\varphi\in X$ and all sufficiently small
$\varepsilon>0$,
\[
\mathcal I_\lambda(u_\lambda+\varepsilon\varphi)
\geq\mathcal I_\lambda(u_\lambda).
\]
Thus Lemma~\ref{lem:identification} applies with $t_\varepsilon=1$,
and $u_\lambda$ is a positive weak solution.

By Theorem~\ref{th5.5}, there is also a nonnegative minimizer
$v_\lambda\in\mathcal N_\lambda^-$. Fix $\varphi\in X$.
Since $G(v_\lambda)>0$, we have
$G(v_\lambda+\varepsilon\varphi)>0$ for sufficiently small
$\varepsilon>0$. Set
\[
t_\varepsilon:=t^-(v_\lambda+\varepsilon\varphi).
\]
Lemma~\ref{lem:ray-continuity} gives
$t_\varepsilon\to t^-(v_\lambda)=1$. Since
$(v_\lambda+\varepsilon\varphi) t_\varepsilon\in\mathcal N_\lambda^-$,
minimality on this component and \eqref{eq:global-ray-max} give
\[
\mathcal I_\lambda\bigl(t_\varepsilon(v_\lambda+\varepsilon\varphi)\bigr)
\geq\mathcal I_\lambda(v_\lambda)
\geq\mathcal I_\lambda(t_\varepsilon v_\lambda).
\]
Hence Lemma~\ref{lem:identification} applies to $v_\lambda$ as well.

Finally,
\[
\mathcal I_\lambda(u_\lambda)<0<\eta_0
\leq\mathcal I_\lambda(v_\lambda),
\]
so the two weak solutions are distinct. Their positivity and the
integrability of the singular term have been established in
Lemma~\ref{lem:identification}. This completes the proof.
\end{proof}

\begin{corollary}\label{cor:small-branch}
The first solution satisfies
\begin{equation}\label{eq:small-rate}
\|u_\lambda\|_X
\leq\left(\frac{4C_f\lambda}{3}\right)^{1/(p^+-s^-)}.
\end{equation}
In particular, $u_\lambda\to0$ strongly in $X$ as $\lambda\downarrow0$.
\end{corollary}

\begin{proof}
Since $\|u_\lambda\|_X<r_0$,
\[
G(u_\lambda)
\leq C_g\|u_\lambda\|_X^{q^-}
\leq\frac14\|u_\lambda\|_X^{p^+}
\leq\frac14\rho(u_\lambda).
\]
The Nehari identity therefore implies
\[
\frac34\|u_\lambda\|_X^{p^+}
\leq\frac34\rho(u_\lambda)
\leq\lambda F(u_\lambda)
\leq\lambda C_f\|u_\lambda\|_X^{s^-}.
\]
Dividing by the positive factor $\|u_\lambda\|_X^{s^-}$ proves
\eqref{eq:small-rate}.
\end{proof}

The preceding result relies on the positivity assumption imposed
on the singular weight $f$. We now consider the more general case
where the weight functions are nonnegative and may vanish
on subsets of $\Omega$. The following proposition shows that the
multiplicity result remains valid in this setting, although
positivity almost everywhere of the solutions is guaranteed only where $f$ is positive.

\begin{samepage}
\begin{proposition}\label{prop:localized}
Assume {\rm(A1),(A2)} and replace {\rm(A3)} by
\[
f,g\in C(\overline\Omega),\qquad
f,g\geq0,\qquad f\not\equiv0,\quad g\not\equiv0.
\]
For every sufficiently small $\lambda>0$, the same construction gives
two distinct nontrivial nonnegative solutions of the equation and
boundary condition in \eqref{P}, with negative and positive energies.
Both solutions are positive almost everywhere on $\{f>0\}$.
Here the singular term is interpreted as
\[
T_u(x):=
\begin{cases}
f(x)u(x)^{-\delta(x)},&f(x)>0,\\
0,&f(x)=0,
\end{cases}
\]
and the weak formulation uses $T_u$ in place of $fu^{-\delta(x)}$.
No assertion of positivity throughout $\Omega$ is made under these
weaker weight assumptions.
\end{proposition}
\end{samepage}

\begin{proof}
The embedding estimates and Lemma~\ref{lem:sphere} are unchanged.
For directions with $F(u)>0$, Lemma~\ref{lem:fibering} applies without
modification. If $F(u)=0<G(u)$, the function
\[
t^{-p^+}J_\lambda(tu)
=\sum_{i=1}^N\int_\Omega a_i(u,x)t^{p_i(x)-p^+}\,dx
-\int_\Omega g_u t^{q(x)-p^+}\,dx
\]
is strictly decreasing from a positive value or $+\infty$ to
$-\infty$. It has exactly one zero, on the descending branch and
beyond $r_0/\|u\|_X$. The corresponding fibering map has a positive
global maximum there. If $F(u)=G(u)=0$, there is no Nehari
intersection on that ray.

The first minimization argument uses a smooth test function with
$F(\psi)>0$, which exists because $f\not\equiv0$.
The second minimization argument and the continuity of the descending
radial parameter remain valid on $\{G>0\}$. Finally, the Fatou
argument in Lemma~\ref{lem:identification} yields positivity on
$\{f>0\}$ and integrability of $T_u|\varphi|$. On $\{f=0\}$ the
singular energy and all its directional increments vanish. Thus the
dominated convergence argument still proves the weak equation.
\end{proof}

\begin{remark}\label{rem:compact-weight}
Proposition~\ref{prop:localized} applies, in particular, when both
weights have compact support in $\Omega$. It separates the existence
argument from the additional problem of propagating positivity beyond
the support of the singular weight. Such propagation cannot be
inferred from an isotropic $p(x)$-Laplacian Harnack inequality.
Strong maximum principles in anisotropic settings depend on the
operator and on its directional structure; see \cite{Vetois2012}.
Likewise, Theorem~\ref{thm:two-solutions} asserts positivity almost
everywhere, not continuity, a pointwise positive representative, or a
positive lower bound on every compact subset. Those stronger
properties require additional regularity arguments.
\end{remark}

\section{Conclusion and Perspectives}\label{sec:conclusion}

We have established the existence of two distinct positive weak
solutions with opposite energy signs for a class of singular
anisotropic elliptic problems with variable exponents, for
sufficiently small $\lambda>0$. A key feature of our approach is
the explicit and uniform construction of the constants
$r_0$, $\lambda_0$, and $\eta_0$ in Lemma~\ref{lem:sphere},
which simultaneously ensures $J_\lambda(u)>0$ and
$\mathcal I_\lambda(u)\geq\eta_0>0$ on the sphere
$\|u\|_X=r_0$. This provides a quantitative energy barrier
and the separation needed for the Nehari geometry despite
the absence of a common homogeneity and the presence of
a singular term. Combined with constrained minimization
and one-sided variational arguments, this yields multiplicity
without requiring global $C^1$ differentiability of the
energy or regularization of the singularity. The extension
to nonnegative weights further distinguishes existence
from positivity throughout the domain.

Several questions naturally arise from these results.
\begin{enumerate}
\item \emph{Optimal parameter range and continuation of solutions.}
The number $\lambda_0$ obtained here is a sufficient smallness
threshold, not an optimal existence threshold. A natural problem is
to study
\[
\Lambda^*:=\sup\{\lambda>0:\text{problem \eqref{P} has a positive weak solution}\},
\]
including whether $\Lambda^*$ is finite, whether a solution exists at
$\Lambda^*$, and how the two variational families behave as the
parameter increases. This project would require a priori bounds and
an analysis of the solution set; the present proof does not establish
continuous or ordered branches.

\item \emph{Small-parameter asymptotics and limiting profiles.}
Estimate \eqref{eq:small-rate} already shows that the negative-energy
solution tends to zero in $X$. It would be useful to determine whether
this rate is sharp and whether an appropriate normalization produces
a nontrivial limiting profile. For the positive-energy solutions, a
related problem is to establish convergence, along suitable sequences
as $\lambda\downarrow0$, to a nontrivial solution of
\[
-\sum_{i=1}^N\partial_i
\bigl(|\partial_i v|^{p_i(x)-2}\partial_i v\bigr)
+\sum_{i=1}^N|v|^{p_i(x)-2}v
=g(x)v^{q(x)-1}.
\]
The absence of a common homogeneous scaling makes the identification
of such profiles a separate question.

\item \emph{Regularity and boundary behavior.}
The present results are formulated in the weak sense. Under additional
regularity and growth assumptions on the exponents, one may study
local boundedness, H\"older continuity, and quantitative lower bounds
on compact subsets. A further project is to determine how the
boundary decay depends on the singular exponent, the directional
growth, and the vanishing rate of $f$ near $\partial\Omega$.
These questions would clarify when positivity almost everywhere can
be strengthened to pointwise positivity of a continuous
representative.

\item \emph{Sharper embeddings and critical growth.}
The proof uses the conservative subcritical scale associated with
$p_m(x)$. It is natural to investigate the larger growth ranges
available from genuinely anisotropic embeddings, whose formulation
involves the harmonic mean of the directional exponents. The critical
case would additionally require a concentration--compactness analysis
and an energy threshold preventing loss of compactness. Extending the
argument to coupled systems would raise related questions concerning
the interaction of several singular components.
\end{enumerate}

These directions extend the existence theory developed here while
keeping separate the questions of multiplicity, positivity,
regularity, and parameter dependence.

\medskip

\noindent {\bf Acknowledgment:} Part of this work was carried out while the second author was visiting KTH Royal Institute of Technology in Stockholm, Sweden. SK gratefully acknowledges
the hospitality received during the visit.

\medskip

\noindent {\bf  AI declaration:} 
The initial draft of this manuscript was prepared without artificial intelligence assistance. During the revision, the authors used ChatGPT (GPT‑6) for proofreading, language editing and suggesting directions for
future research. The tool suggested items (2) and (4) in Section \ref{sec:conclusion}, and also some relevant references like \cite{ChemsEddineRagusaRepovs2024,ChemsEddineRagusa2025,ChemsEddineRagusa2026} which were previously unknown to the authors.

\medskip

\noindent {\bf  Data availability statement:} No datasets were generated or analyzed during this study.

\medskip
\noindent {\bf  Funding and/or Conflicts of interests/Competing interests:} The authors declare that there are no financial, competing or conflict of interests.

\end{document}